\documentstyle[12pt]{article}
\def\mymedskip{\vskip\medskipamount}
\def\mymedbreak{\par \ifdim\lastskip<\medskipamount
  \removelastskip \penalty-100 \mymedskip \fi}
\def\myaftermedspace{\par \ifdim\lastskip<\medskipamount
  \removelastskip \penalty55\mymedskip\fi}
\newcommand{\eop}{{\unskip\nobreak\hfil\penalty50
          \hskip2em\hbox{}\nobreak\hfil$\Box$
          \parfillskip=0pt \finalhyphendemerits=0 \par}}
\newenvironment{proof}%
{\mymedbreak{\noindent\bf Proof:\enspace}}{\eop\myaftermedspace}
{\mymedbreak{\noindent\bf Proof of Theorem #1:\enspace}}{\eop\myaftermedspace}
\newtheorem{teor}{Theorem}[section]
\newtheorem{defi}[teor]{Definition}
\newtheorem{examp}[teor]{Example}
\newtheorem{remark}[teor]{Remark}

\newtheorem{con}[teor]{Conjecture}

\newcommand{\beq}{\begin{equation}}
\newcommand{\eeq}{\end{equation}}
\newcommand{\beql}[1]{\begin{equation} \label{#1}}
\newcommand{\eeql}{\end{equation}}
\newcommand{\beqa}{\begin{eqnarray*}}
\newcommand{\eeqa}{\end{eqnarray*}}
\newcommand{\beqal}[1]{\begin{eqnarray} \label{#1}}
\newcommand{\eeqal}{\end{eqnarray}}
\newcommand{\beqan}{\begin{eqnarray}}
\newcommand{\eeqan}{\end{eqnarray}}
\newcommand{\bpf}{\begin{proof}}
\newcommand{\epf}{\end{proof}}

\newcommand{\cE}{{\cal E}}

\newcommand{\cO}{{\cal O}}
\newcommand{\cR}{{\cal R}}
\newcommand{\cH}{{\cal H}}

\newcommand{\bF}{{\bf F}}
\newcommand{\bT}{{\bf T}}

\newcommand{\PG}{{\rm PG}}
\newcommand{\PGL}{{\rm PGL}}
\newcommand{\PgammaL}{{\rm P\Gamma L}}

\newcommand{\Tr}{{\rm Tr}}

\begin{document}
\begin{titlepage}
\title{Pseudocyclic association schemes arising from the actions of $\PGL(2,2^m)$ and $\PgammaL(2,2^m)$}
\date{In Memory of Jack van Lint}
\author{%
Henk D. L.\ Hollmann\\Philips Research Laboratories\\
Prof. Holstlaan 4, 5656 AA Eindhoven\\The Netherlands\\email: {\tt
henk.d.l.hollmann@philips.com}\\
\\Qing Xiang\\Department of Mathematical
Sciences\\University of Delaware\\Newark, DE 19716, USA\\email: {\tt
xiang@math.udel.edu}%
}
\maketitle
\begin{abstract}
The action of $\PGL(2,2^m)$ on the set of exterior lines to a nonsingular conic in $\PG(2,2^m)$ affords an association scheme, which was shown to be pseudocyclic in \cite{henkthesis}. It was further conjectured in \cite{henkthesis} that the orbital scheme of $\PgammaL(2,2^m)$ on the set of exterior lines to a nonsingular conic in $\PG(2,2^m)$ is also pseudocyclic if $m$ is an odd prime. We confirm this conjecture in this paper. As a by-product, we obtain a class of Latin square type strongly regular graphs on nonprime-power number of points.
\end{abstract}
\end{titlepage}
\newpage

\section{\label{intro}Introduction}

Let $X$ be a finite set. A (symmetric) {\it association scheme} with $d$ classes on $X$ is a partition of
$X\times X$ into sets $R_0$, $R_1, \ldots , R_d$ (called {\it associate classes} or {\it relations}) such that
\begin{enumerate}
\item $R_0=\{(x,x) \mid  x\in X\}$ (the diagonal relation),
\item $R_i$ is symmetric for $i=1,2,\ldots ,d$,
\item for all $i,j,k$ in $\{0,1,2,\ldots ,d\}$ there is an integer $p_{ij}^k$ such that, for all $(x,y)\in R_k$,
$$|\{z\in X \mid (x,z)\in R_i\; {\rm and}\; (z,y)\in R_j\}|=p_{ij}^k.$$
\end{enumerate}
We denote such an association scheme by $(X, \{R_i\}_{0\leq i\leq d})$. Elements $x$ and $y$ of $X$ are called {\it $i$-th associates} if $(x,y)\in R_i$. The numbers $p_{ij}^k$, $0\leq k,i,j\leq d$, are called the {\it intersection parameters} of the scheme. That $p_{ii}^0$ exists means that there is a constant number of $i$-th associates of any element of $X$, which is usually denoted by $n_i$. The numbers $n_0,n_1,\ldots ,n_d$ are called the {\it valencies} (or {\it degrees}) of the scheme. We have
\begin{enumerate}
\item $n_0=1$, $n_0+n_1+\cdots +n_d=|X|,$
\item $p_{0j}^k=\delta_{j,k}$ (Kronecker delta), $p_{ij}^0=\delta_{i,j}n_j$,
\item $p_{ij}^k=p_{ji}^k$, $p_{ij}^kn_k=p_{ik}^jn_j$.
\end{enumerate}
For $i\in \{0,1,\ldots ,d\}$, let $A_i$ be the adjacency matrix of the relation $R_i$, that is, the rows and columns of $A_i$ are both indexed by $X$ and
$$(A_i)_{xy}:=\biggm\{
\begin{array}{c} 1 \quad \mbox{ if }\quad (x,y)\in R_i, \\
                 0 \quad \mbox{ if }\quad (x,y)\notin R_i. \\
\end{array} $$
The matrices $A_i$ are symmetric $(0,1)$-matrices and
$$A_0=I, \; A_0+A_1+\cdots +A_d=J,$$
where $J$ is the all one matrix of size $|X|$ by $|X|$.

By the definition of an association scheme, we have
$$A_iA_j=\sum_{k=0}^d p_{ij}^kA_k $$
for any $i,j\in \{0,1,\ldots ,d\}$. So $A_0,A_1,\cdots , A_d$ form a basis of the commutative algebra generated
by $A_0,A_1,\cdots , A_d$ over the reals (which is called the {\it Bose-Mesner algebra} of the association
scheme). Moreover this algebra has a unique basis $E_0,E_1,\cdots , E_d$ of primitive idempotents; one of the
primitive idempotents is $\frac {1}{|X|}J$.  So we may assume that $E_0=\frac {1} {|X|}J$. Let $m_i={\rm
rank}\;E_i$. Then
$$m_0=1,\; m_0+m_1+\cdots +m_d=|X|.$$
The numbers $m_0,m_1,\ldots ,m_d$ are called the {\it multiplicities} of the scheme. Since we have two bases of
the Bose-Mesner algebra, we may consider the transition matrices between them. Define
$P=\left(p_j(i)\right)_{0\le i,j\le d}$ (the {\it first eigenmatrix}) and $Q=\left(q_j(i)\right)_{0\le i,j\le d}$
(the {\it second eigenmatrix}) as the $(d+1)\times (d+1)$ matrices with rows and columns indexed by $0,1,2,\ldots
,d$ such that
$$(A_0,A_1, \ldots ,A_d)=(E_0,E_1, \ldots ,E_d)P,$$
and
$$|X|(E_0,E_1, \ldots ,E_d)=(A_0,A_1, \ldots ,A_d)Q.$$
Of course, we have
$$P=|X|Q^{-1}, \;\; Q=|X|P^{-1}.$$
Note that $\{p_j(i)\ |\ 0\le i\le d\}$ is the set of eigenvalues of $A_j$ and the zeroth row and column of $P$ and $Q$ are as indicated below.
$$P=\left(\matrix{
1&n_1&\cdots&n_d\cr
1\cr
\vdots & & \cr
1
}\right),\;\; Q=\left(\matrix{
1&m_1&\cdots&m_d\cr
1\cr
\vdots & & \cr
1
}\right)$$

Before we proceed further, we give some examples of association schemes.

\begin{examp}\label{tranexamp}
Let $X$ be a finite set and let $G$ be a group acting transitively on $X$. We say that $G$ acts {\em generously transitively} on $X$ if the orbits of the induced action of $G$ on $X\times X$ are all symmetric. (The orbits of $G$ on $X\times X$ are usually called the {\em orbitals} of the action of $G$ on $X$.) It is clear that if $G$ acts generously transitively on $X$, then the orbitals of $G$ on $X$ can be taken as the relations of an association scheme, which will be called the {\em orbital scheme} of $G$ on $X$. The next example arises in this way.
\end{examp}

\begin{examp}\label{cycexamp}
We consider {\em cyclotomic schemes} defined as follows. Let $q$ be a prime power and let $q-1=ef$ with $e>1$. Let $C_0$ be the subgroup of the multiplicative group of $\bF_q$ of index $e$, and let $C_0,C_1,\ldots ,C_{e-1}$ be the cosets of $C_0$. We require $-1\in C_0$. Define $R_0=\{(x,x) : x\in \bF_q\}$, and for $i\in \{1,2,\ldots ,e\}$, define $R_i=\{(x,y)\mid x,y\in \bF_q, x-y\in C_{i-1}\}$. Then $(\bF_q, \{R_i\}_{0\leq i\leq e})$ is an $e$-class symmetric association scheme (the $R_i$ are the orbitals of the action of $G$ on $\bF_q$, where $G=\{x\mapsto ax+b\mid a\in C_0, b\in \bF_q\}$). The intersection parameters of the cyclotomic scheme are related to the cyclotomic numbers (\cite[p.~25]{st}). Namely, for $i,j,k\in \{1,2,\ldots ,e\}$, given $(x,y)\in R_k$,
\begin{equation}\label{cycparam}
p_{ij}^k=|\{z\in \bF_q\mid x-z\in C_{i-1}, y-z\in C_{j-1}\}|=|\{z\in C_{i-k}\mid 1+z\in C_{j-k}\}|.
\end{equation}
The first eigenmatrix $P$ of this scheme is the following $(e+1)$ by $(e+1)$ matrix (with the rows of $P$ arranged in a certain way)
$$P=\left(\matrix{
1&f&\cdots&f\cr
1\cr
\vdots & & P_0\cr
1
}\right)$$
with $P_0=\sum_{i=1}^{e}\eta_iC^i$, where $C$ is the $e$ by $e$ matrix:
$$C=\left(\matrix{
& 1\cr & & 1\cr & & & \ddots\cr & & & & 1\cr 1}\right)$$ and $\eta_i=\sum_{\beta\in C_i}\psi(\beta)$, $1\leq
i\leq e$, for a fixed nontrivial additive character $\psi$ of $\bF_q$. See \cite{bandm} for more details.
\end{examp}

Next we introduce the notion of a pseudocyclic association scheme.

\begin{defi}
Let $(X, \{R_i\}_{0\leq i\leq d})$ be an association scheme. We say that $(X, \{R_i\}_{0\leq i\leq d})$ is {\it pseudocyclic} if there exists an integer $t$ such that $m_i=t$ for all $i\in \{1,\cdots , d\}$.
\end{defi}

The following theorem gives combinatorial characterizations for an association scheme to be pseudocyclic.

\begin{teor}\label{pseudocyc}
Let $(X, \{R_i\}_{0\leq i\leq d})$ be an association scheme, and for $x\in X$ and $1\leq i\leq d$, let $R_i(x)=\{y\mid (x,y)\in R_i\}$. Then the following are equivalent.\\
(1). $(X, \{R_i\}_{0\leq i\leq d})$ is pseudocyclic.\\
(2). For some constant $t$, we have $n_j=t$ and $\sum_{k=1}^{d}p_{kj}^k=t-1$, for $1\leq j\leq d$.\\
(3). $(X, {\mathcal B})$ is a $2-(v,t,t-1)$ design, where ${\mathcal B}=\{R_i(x)\mid x\in X, 1\leq i\leq d\}$.
\end{teor}

For a proof of this theorem, we refer the reader to \cite[p.~48]{bcn} and \cite[p.~84]{henkthesis}. Part (2) in the above theorem is very useful. For example, we may use it to prove that the cyclotomic scheme in Example~\ref{cycexamp} is pseudocyclic. The proof goes as follows. First, the nontrivial valencies of the cyclotomic scheme in Example~\ref{cycexamp} are all equal to $f$. Second, by (\ref{cycparam}) and noting that $-1\in C_0$, we have
\begin{eqnarray*}
\sum_{k=1}^ep_{kj}^k&=&\sum_{k=1}^e|\{z\in C_{0}\mid 1+z\in C_{j-k}\}|\\
&=&|C_0|-1=f-1\\
\end{eqnarray*}

Pseudocyclic schemes can be used to construct strongly regular graphs and distance regular graphs of diameter 3
(\cite{bm}, \cite[p.~388]{bcn}). In view of this, it is of interest to construct pseudocyclic association
schemes, as remarked by the authors of \cite{bcn} (see \cite[p.~389]{bcn}). The cyclotomic schemes are examples of
pseudocyclic association schemes on prime-power number of points. Very few examples of pseudocyclic association
schemes on nonprime-power number of points are currently  known (see \cite{mathon}, \cite[p.~390]{bcn} and
\cite{henkthesis}). One class of such examples comes from the action of $\PGL(2,2^m)$ on the set of exterior
lines to a nonsingular conic in $\PG(2,2^m)$. We will give a quick review of this class of association schemes in
Section 2, and also include a proof of the pseudocyclicity of these association schemes. In \cite{henkthesis}, it
was further conjectured that the orbital scheme of $\PgammaL(2,2^m)$ on the set of exterior lines to a
nonsingular conic in $\PG(2,2^m)$ is also pseudocyclic if $m$ is an odd prime. We will confirm this conjecture in
Section 3. As a by-product, we obtain a class of Latin square type strong regular graphs on nonprime-power number
of points.

\section{The Elliptic Schemes}

In the rest of this paper, we always assume that $q=2^m$, where $m$ is a positive integer. Let
$$\cO=\{(\xi, \xi^2, 1)\mid \xi\in \bF_q\}\cup\{(0,1,0)\}.$$
Then $\cO$ is a nonsingular conic in $\PG(2,q)$. A line of $\PG(2,q)$ is called {\it exterior} (resp. {\it
secant}) if it meets $\cO$ in 0 (resp. 2) points. Let $\cE$ (resp. $\cH$) be the set of exterior (resp. secant)
lines to $\cO$. Then $$|\cE|=\frac {q(q-1)} {2},\; \mbox{and}\; |\cH|=\frac {(q+1)q} {2}.$$ The subgroup of
$\PGL(3,q)$ fixing $\cO$ setwise is isomorphic to $\PGL(2,q)$ (cf. \cite[p.~158]{hirsch}). Hence $\PGL(2,q)$ acts
on $\cE$ and $\cH$, respectively. Moreover, it is shown in \cite{hxmay2004} that $\PGL(2,q)$ acts generously
transitively on both $\cE$ and $\cH$. Therefore we obtain two association schemes, one on $\cE$ and the other on
$\cH$. The association scheme on $\cE$ will be called the {\it elliptic} scheme, and the association scheme on
$\cH$ is called the {\it hyperbolic} scheme.

Since the point $(1,0,0)$ is the nucleus of $\cO$ (i.e., the point at which all tangent lines to $\cO$ meet), we
see that each line in $\cE\cup\cH$ can be written as $(1,x,y)^{\perp}=\{(a_0,a_1,a_2)\in \bF_q^3\mid
a_0+a_1x+a_2y=0\}$ for some $x,y\in\bF_q$. Let $\Tr: \bF_q\rightarrow \bF_2$ be the trace map. Also for $e\in\bF_2$
we define
$$\bT_e=\{x\in\bF_q\mid \Tr(x)=e\},$$
and $\bT_e^*=\bT_e\setminus\{0\}$. Then $(1,x,y)^{\perp}$ is in $\cE$ (resp. $\cH$) if and only if $\Tr(xy)=1$ (resp. $\Tr(xy)=0$). Given two lines $\ell=(1,x,y)^{\perp}$ and $m=(1,z,u)^{\perp}$, we define
$${\hat \rho}(\ell,m)=x^2u^2+y^2z^2 +(x+z)(y+u).$$
We remark that the function ${\hat \rho}$ comes from the cross-ratio of four points on a projective line (see
\cite{hxmay2004} for details). The following theorem in \cite{hxmay2004} gives a simple description of the
orbitals of the action of $\PGL(2,q)$ on $\cE$ by using the function ${\hat \rho}$.

\begin{teor}\label{description}
The orbitals of the action of $\PGL(2,q)$ on $\cE$ are $\Gamma_0$ (the diagonal class), and $\Gamma_{a}=\{(\ell, m)\mid {\hat \rho}(\ell, m)=a\}$ for all $a\in \bT_0^*$.\end{teor}

There is a similar description of the orbitals of $\PGL(2,q)$ on $\cH$ in \cite{hxmay2004}. Since we are only concerned with the elliptic scheme in this paper, we omit that description.

The pair $(\cE, \{\Gamma_a\})$ is an association scheme on $\cE$ with $\frac {(q-2)}{2}$ classes. The intersection parameters of this scheme are computed in \cite{hxmay2004}. For $a,b,c\in \bT_0^*$, given $(\ell, m)\in\Gamma_c$, we use $p_{a,b}^c$ to denote $|\{k\in \cE \mid (\ell,k)\in \Gamma_a\; {\rm and}\; (k,m)\in \Gamma_c\}|$. We have:

\begin{teor}\label{parameters}
Let $a,b,c\in \bT_0^*$. Then for any $v\in \bT_1$,
\beql{pexp} p^c_{a,b}=\left\{ \begin{array}{ll}
                                1+2\delta_{\Tr(ac),1}, \; & \mbox{if $a+b+c=0$;} \\
                                \sum_\tau |\{z\in\bF_q \mid z^2+z=v+ac/\tau^2\}|, & \mbox{otherwise,}
                            \end{array}
                    \right.
\eeql
where the last sum is over the two elements $\tau\in \bF_q$ satisfying $\tau^2+\tau=a+b+c$. Also for all $a\in \bT_0^*$, the valency $n_a=q+1$.
\end{teor}

The association scheme $(\cE, \{\Gamma_a\})$ is pseudocyclic. This is already known in \cite{henkthesis}. For convenience of the reader, we include a proof here.

\begin{teor}\label{ellipseudo}
The association scheme $(\cE, \{\Gamma_a\})$ is pseudocyclic.
\end{teor}

\begin{proof}
By Theorem~\ref{parameters}, the nontrivial valencies of the association scheme $(\cE, \{\Gamma_a\})$ are all
equal to $q+1$. By Part (2) of Theorem~\ref{pseudocyc}, it suffices to prove that $\sum_{a\in
\bT_0^*}p_{a,b}^a=q$ for all $b\in \bT_0^*$.

By Theorem~\ref{parameters}, for $a,b\in \bT_0^*$, we have
$$p_{a,b}^a=\sum_{\tau^2+\tau=b}(1-(-1)^{\Tr(a/\tau)}).$$
Fixing $\tau\in\bF_q\setminus\{0,1\}$ with $\tau^2+\tau=b$, we have
\beqa \sum_{a\in \bT_0^*} p^a_{a,b} &=& \sum_{a\in \bT_0^*}(1-(-1)^{\Tr(a/\tau)}+1-(-1)^{\Tr(a/(\tau+1))}) \\
  &=& 2(q/2 -1)- \sum_{a\in \bT_0^*}((-1)^{\Tr(a/\tau)}+(-1)^{\Tr(a/(\tau+1))})\\
   &=& 2(q/2 -1)-(-1-1)\\
   &=& q
\eeqa
This completes the proof.
\end{proof}

\section{Pseudocyclic fusion schemes of the elliptic schemes}

As we have seen in the last section, the elliptic scheme $(\cE,\{\Gamma_a\})$ is pseudocyclic. In this section, we will consider the fusion scheme of $(\cE,\{\Gamma_a\})$ obtained by merging the classes $\Gamma_a$ via the Frobenius automorphism $x\mapsto x^2$ of $\bF_q$. Specifically, for $a\in \bT_0^*$, define
$$\Delta_a=\cup_{i\in C_a} \Gamma_i,$$
where $C_a:=\{a, a^2, a^4, \ldots, a^{2^{m-1}}\}$. Let $\cR$ be a set of orbit representatives of $\bT_0^*$ under
the action of the Frobenius automorphism. Then $\Delta_0:=\Gamma_0$, and $\Delta_a$, $a\in\cR$ are the orbitals
of $\PgammaL(2,q)$ on $\cE$. Therefore $(\cE, \{\Delta_a\})$ is also an association scheme. The (nontrivial)
intersection parameters of this fusion scheme will be denoted by $P^c_{a,b}$, where $a,b,c\in \cR$. We have for
$a,b,c\in\cR$,
\[ P^c_{a,b}=\sum_{e\in C_a} \sum_{f\in C_b} p^g_{e,f},\]
where $g\in C_c$. (This is independent of the choice of $g\in C_c$.)

Now, if $m$ is prime, then each $C_a$, $a\in \cR$, has size $m$, so the nontrivial valencies of the fusion scheme $(\cE, \{\Delta_a\})$ are all equal to $m(q+1)$. Hollmann \cite[p.~133]{henkthesis} made the following conjecture.

\begin{con}\label{pseudoconj}
If $m$ is an odd prime, then $(\cE, \{\Delta_a\})$ is pseudocyclic.
\end{con}

As far as we know, there is no published proof of this conjecture. There is one sentence on Page 390 of \cite{bcn} stating the above conjecture as a fact. But this was not backed up by a proof.

Note that the nontrivial valencies of $(\cE, \{\Delta_a\})$ are all equal to $m(q+1)$ when $m$ is prime. So in
order to prove Conjecture~\ref{pseudoconj}, by Part (2) of Theorem~\ref{pseudocyc}, we need to show that
\beql{original} \sum_{c\in \cR}P_{c,c}^b=m(q+1)-1, \eeql for all $b\in\cR$. (Here we implicitly used the fact
that $P_{c,c}^b=P_{c,b}^c$ since all nontrivial valencies are all equal when $m$ is prime.) Simplifying the left
hand side of (\ref{original}), we see that (\ref{original}) is equivalent to \beql{pc1}
\sum_{k=0}^{m-1}\sum_{c\in \bT_0^*} p^b_{c,c^{2^k}}=m(q+1)-1. \eeql Now, the $k=0$ term of the left hand side of
(\ref{pc1}) is equal to $q$ as seen in the proof of Theorem~\ref{ellipseudo}. So in order to prove (\ref{pc1}),
we have to show that \beql{pc2}\sum_{k=1}^{m-1} \sum_{c\in \bT_0^*} p^b_{c,c^{2^k}}=(m-1)(q+1), \eeql for all
$b\in\bT_0^*$.

We will prove a stronger result:
\begin{teor}\label{strong}
Let $m$ be an odd integer, and let $k$ be any integer in $\{1,2,\ldots ,m-1\}$ satisfying $\gcd(k,m)=1$. Write $\sigma=2^k$. Then
\beql{pc0} \sum_{c\in \bT_0^*} p^b_{c,c^\sigma} =q+1, \eeql
for all $b\in\bT_0^*$.
\end{teor}

The most important ingredient in our proof of Theorem~\ref{strong} is a family of polynomials $H_{\alpha, \gamma}(X)$ introduced in \cite{hxpermpoly}. In fact we discovered these polynomials while working on a proof of Theorem~\ref{strong}. We now define the polynomials $H_{\alpha, \gamma}(X)$ and quote the main theorem from \cite{hxpermpoly}.

Let $m\geq1$ be an integer, let $k$ be any integer in $\{1,\ldots,
m-1\}$ with $\gcd(k,m)=1$, and let $r\in\{1,\ldots, m-1\}$ be such
that $kr\equiv 1$ (mod $m$). Write $\sigma=2^k$ and use $\Tr(X)$ to denote the following polynomial in $\bF_2[X]$.
$$\Tr(X):=X+X^2+\cdots+X^{2^{m-1}}.$$
For $\alpha, \gamma$ in $\{0,1\}$, we define the polynomial
\[H_{\alpha,\gamma}(X):= \gamma \Tr(X) + \frac{\left(\alpha \Tr(X) + \sum_{i=0}^{r-1} X^{\sigma^i}\right)^{\sigma+1}} {X^2}.\]
(Note that $H_{\alpha,\gamma}(X)$ is indeed a polynomial in $X$ with coefficients in $\bF_2$ and $H_{\alpha,\gamma}(0)=0$. Also see \cite{hxpermpoly} for connections between $H_{\alpha,\gamma}(X)$ and the Dickson polynomials.)

The following is the main theorem from \cite{hxpermpoly}.

\begin{teor}\label{mainthm} Let $m, k$ be positive integers with $\gcd(k,m)=1$,
let $r\in\{1,\ldots, m-1\}$ be such that $kr\equiv 1 \;(\bmod
\;m)$, and let $\alpha, \gamma\in \{0,1\}$. Then the mapping $H_{\alpha, \gamma}: x\mapsto H_{\alpha,\gamma}(x)$, $x\in\bF_q$, maps $\bT_0$ bijectively to $\bT_0$, and maps $\bT_1$ bijectively to $\bT_{r+(\alpha+\gamma)m}$. In particular, the polynomial $H_{\alpha,\gamma}(X)$ is a permutation polynomial of $\bF_q$ if and only if $r+(\alpha+\gamma)m \equiv 1$ {\em (mod 2)}.
\end{teor}

We are now ready to give the proof of Theorem~\ref{strong}.

\vspace{0.1in}

\noindent{\bf Proof of Theorem~\ref{strong}:} Recall that from Theorem~\ref{parameters}, for $b,c\in \bT_0^*$,
\begin{eqnarray*}
p^b_{c,c^\sigma}=\left\{ \begin{array}{ll}
                                1+2\delta_{\Tr(bc),1}, \; & \mbox{if $c^\sigma+c+b=0$;} \\
                                \sum_{\tau^2+\tau=c^\sigma+c+b} |\{z\in\bF_q \mid z^2+z=v+bc/\tau^2\}|, & \mbox{if $c^\sigma+c+b\neq 0$,}
                            \end{array}
                    \right.
\end{eqnarray*}
where $v$ is any element with $\Tr(v)=1$. Since $b\in\bT_0^*$ and $m$ is odd, we can find a unique $c_0\in \bT_0^*$ such that $c_0^\sigma +c_0=b$. So
\begin{eqnarray*}
\sum_{c\in \bT_0^*} p^b_{c,c^\sigma} &=& 1+2\delta_{\Tr(bc_0),1}+2\sum_{c\in \bT_0^*,\; c^\sigma +c+b\neq 0}\sum_{\tau^2+\tau=c^\sigma+c+b}\delta_{\Tr(bc/\tau^2), 1}\\
&=&1+2|\{(c,\tau)\in \bF_q^*\times \bF_q^*\mid \tau^2+\tau=c^\sigma+c+b, \Tr(c)=0, \Tr(bc/\tau^2)=1\}|.
\end{eqnarray*}
For convenience, we define
$$N_{k}(b):=|\{(c,\tau)\in \bF_q^*\times \bF_q^*\mid \tau^2+\tau=c^\sigma+c+b, \Tr(c)=0, \Tr(bc/\tau^2)=1\}|.$$
Our goal is to prove that $N_{k}(b)=q/2$ for all $b\in\bT^*_0$.

For later use, we define the polynomial
\[ f(X):= \sum_{i=0}^{r-1} X^{\sigma^i}\in \bF_2[X], \]
where $r$ is an integer satisfying $kr\equiv 1$ (mod $m$).

Since $b\in \bT_0^*$ and $m$ is odd, we can write $b=\beta +\beta^2$ with $\beta\in\bT_0^*$. Then the equation
$\tau^2+\tau=c^\sigma +c+b$ involved in the definition of $N_{k}(b)$ becomes \beql{newequ} c^\sigma +c=(\beta
+\tau)+(\beta +\tau)^2. \eeql Noting that $m$ is odd, we see that for any $\tau\in\bF_q$, there is a unique
solution $c\in\bT_0$ of (\ref{newequ}), namely
 \[c=f(\tau +\beta)+r\Tr(\tau +\beta)=f(\tau +\beta)+r\Tr(\tau),\]
where in the last equality we used the fact that $\beta\in \bT_0$. Therefore we have
\begin{eqnarray*}
N_{k}(b)=\left\{ \begin{array}{ll}
                                |\{\tau\in\bF_q^* :  \frac{b (f(\tau +\beta)+\Tr(\tau))}{\tau^2}\in \bT_1\}|, & \mbox{if $r$ is odd;} \\
                              |\{\tau\in\bF_q^* : \frac{b f(\tau +\beta)}{\tau^2}\in \bT_1\}|  , & \mbox{if $r$ is even.}
                            \end{array}
                    \right.
\end{eqnarray*}
We will consider the $r$ odd case and the $r$ even case separately.

\noindent{\bf Case 1}. $r$ is odd. Let $x=b/\tau^2$, where $b=\beta +\beta^2\in \bT_0^*$ and $\tau\in\bF_q^*$. Then
\begin{eqnarray*}
\Tr\left(\frac{b (f(\tau +\beta)+\Tr(\tau))}{\tau^2}\right)&=&\Tr\left(x\sum_{i=0}^{r-1}(\beta+\sqrt{b/x})^{\sigma^i}+x\Tr(b/x)\right) \\
&=&\Tr\left(\sum_{i=0}^{r-1} x^2(\beta^2+b/x)^{\sigma^i}\right)+\Tr(x)\Tr(b/x) \\
&=&\Tr\left(\sum_{i=0}^{r-1} x^{\sigma^{r-i}}(\beta^2+b/x)\right)+\Tr(x)\Tr(b/x)\\
&=&\Tr\left((\beta^2+b/x)(f(x)+x^2+x)\right)+\Tr(x)\Tr(b/x) \\
&=&\Tr\left(\beta^2(f(x)+\frac {f(x)}{x}+\frac{f(x)^2}{x^2})\right)+\Tr(x)\Tr\left(\frac{b}{x}\right),
\end{eqnarray*}
where in the last step, we used $b=\beta+\beta^2$. Now noting that for $x\in \bF_q^*$,
$$H_{0,0}(x)=f(x)+\frac {f(x)}{x}+\frac{f(x)^2}{x^2}.$$
(One can prove this directly, or see Lemma 3.1 in \cite{hxpermpoly}.) Therefore, in this case, we have
\beql{twoterms}N_{k}(b)=|\{x\in\bT_0^*\mid \beta^2 H_{0,0}(x)\in \bT_1\}|+|\{x\in\bT_1\mid \beta^2 H_{0,0}(x)+b/x\in \bT_1\}|.\eeql

For the first summand in (\ref{twoterms}), noting that $H_{0,0}(0)=0$ and $H_{0,0}$ maps $\bT_0$ to $\bT_0$ bijectively (Theorem~\ref{mainthm}), we have
\begin{eqnarray*}
|\{x\in\bT_0^*\mid \beta^2 H_{0,0}(x)\in \bT_1\}| &=& |\beta^2 \bT_0^*\cap\bT_1|\\
&=& (q/2 -1)-|\beta^2\bT_0^*\cap \bT_0^*|.
\end{eqnarray*}
Since $\bT_0^*$ is a $(q-1,q/2 -1,q/4 -1)$ (Singer) difference set in the cyclic group $\bF_q^*$, and $\beta\neq 0, 1$, we see that $|\beta^2\bT_0^*\cap \bT_0^*|=q/4 -1$. Hence
\[|\{x\in\bT_0^*\mid \beta^2 H_{0,0}(x)\in \bT_1\}|=(q/2-1)-(q/4 -1)=q/4.\]

For the second summand in (\ref{twoterms}), using $b=\beta+\beta^2$, we see that
\[\Tr(\beta^2 H_{0,0}(x)+ b/x)=\Tr(\beta^2(H_{0,0}(x)+1/x +1/x^2)).\]
For any $x\in\bT_1$, we have
\beqa
H_{1,0}(x) &=& \frac {(1+ f(x))^{\sigma +1}} {x^2} \nonumber\\
          &=& 1+f(x)+ (1+f(x))/x +(1 +f(x))^2/x^2 \nonumber \\
          &=& 1+ 1/x +1/x^2 + H_{0,0}(x).
\eeqa
Also by Theorem~\ref{mainthm}, $H_{1,0}$ maps $\bT_1$ bijectively to $\bT_{r+m}=\bT_0$. Hence
\begin{eqnarray*}
|\{x\in\bT_1\mid \beta^2 H_{0,0}(x)+b/x\in \bT_1\}|&=&|\{x\in\bT_1\mid \beta^2 (H_{0,0}(x)+1/x+1/x^2)\in \bT_1\}|\nonumber\\
&=&|\beta^2 \bT_1\cap \bT_1|\nonumber\\
&=&q/4.
\end{eqnarray*}
Therefore we have $N_{k}(b)=q/4 +q/4=q/2.$

\noindent{\bf Case 2.} $r$ is even. This case is similar to Case 1 and actually easier. Let $x=b/\tau^2$. By the
same computations as those in the $r$ odd case, we find that
\[\Tr(\frac{b f(\tau +\beta)}{\tau^2})=\Tr\left(\beta^2H_{0,0}(x)\right).\]
By Theorem~\ref{mainthm}, $H_{0,0}$ maps $\bT_0^*$ bijectively to $\bT_0^*$, and maps $\bT_1$ bijectively to $\bT_r=\bT_0$. Therefore,
\beqa
|\{\tau\in\bF_q^*\mid \frac {b f(\tau +\beta)}{\tau^2}\in \bT_1\}|&=&|\{x\in \bF_q^*\mid \beta^2 H_{0,0}(x)\in\bT_1\}|\nonumber\\
&=&|\{x\in \bT_0^*\mid \beta^2 H_{0,0}(x)\in\bT_1\}|+|\{x\in \bT_1\mid \beta^2 H_{0,0}(x)\in\bT_1\}| \nonumber\\
&=&|\beta^2\bT_0^*\cap \bT_1|+|\beta^2\bT_0\cap \bT_1|\nonumber\\
&=&2|\beta^2\bT_0^*\cap \bT_1|\nonumber\\
&=&q/2.
\eeqa

In summary, in both cases, we have shown that $N_{k}(b)=q/2$ for all $b\in\bT_0^*$. The proof is complete.

\begin{remark} More general results can be proved in the same fashion as above. Let $e,f\in\bF_2$. Define
$$N_{k,e,f}(b):=|\{(c,\tau)\in \bF_q^*\times \bF_q^*\mid \tau^2+\tau=c^\sigma+c+b, \Tr(c)=e, \Tr(bc/\tau^2)=f\}|.$$
Then using the same arguments as those in the proof of Theorem~\ref{strong}, we find that $N_{k,0,0}(b)=q/2 -3$, $N_{k,1,0}(b)=q/2 -1$, and $N_{k,1,1}(b)=q/2$, for all $b\in \bT_0^*$.
\end{remark}

Now we can finish the proof of Conjecture~\ref{pseudoconj}.

\begin{teor}
If $m$ is an odd prime, then $(\cE, \{\Delta_a\})$ is pseudocyclic.
\end{teor}

\begin{proof} Since $m$ is prime, the nontrivial valencies of the scheme are all equal to $m(q+1)$. To finish the proof, we need to prove (\ref{original}) for all $b\in \cR$. As we have seen in the analysis before the statement of Theorem~\ref{strong}, (\ref{original}) is equivalent to (\ref{pc2}). Since $m$ is an odd prime, any integer $k\in\{1,2,\ldots ,m-1\}$ is relatively prime to $m$. So we can apply Theorem~\ref{strong} to obtain
\[\sum_{c\in \bT_0^*} p^b_{c,c^\sigma} =q+1,\]
for all $b\in\bT_0^*$. Now (\ref{pc2}) follows. This completes the proof.
\end{proof}

\section{Latin square type strongly regular graphs}

A {\em strongly regular graph srg} $(v,k,\lambda,\mu)$ is a graph with $v$ vertices that is regular of valency $k$ and that has the following properties:
\begin{enumerate}
\item For any two adjacent vertices $x,y$, there are exactly $\lambda$ vertices adjacent to both $x$ and $y$.
\item For any two nonadjacent vertices $x,y$, there are exactly $\mu$ vertices adjacent to both $x$ and $y$.
\end{enumerate}
It is well known \cite[p.~407]{vanlint} that strongly regular graphs are equivalent to two-class association schemes. An srg $(v,k,\lambda,\mu)$ is said to be of {\it Latin square type} if
$$(v,k,\lambda, \mu)=(n^2, t(n-1), n+t^2-3t, t^2-t),$$
where $1\leq t\leq n+1$. Any Latin square of order $n$ gives rise to a Latin square type srg (actually called Latin square graph in this case) with parameters $(n^2, 3(n-1), n-2, 6)$ (see \cite[p.~273]{vanlint}). Many examples of Latin square type srg on prime-power number of points are known \cite{ma}. In contrast, not too many examples of Latin square type srg on nonprime-power number of points are known.

In \cite{bm}, it was shown that pseudocyclic association schemes can give rise to Latin square type srg. We quote the following theorem from \cite{bm}. A proof can be found in \cite{tf}.

\begin{teor}\label{srg}
Let $(X, \{R_i\}_{0\leq i\leq d})$ be a pseudocyclic association scheme on $dt+1$ points. Then the graph $G$ whose vertex set is $X\times X$, and where two distinct vertices $(x,y)$ and $(x',y')$ are adjacent if and only if $(x,x')\in R_i$ and $(y,y')\in R_i$ for some $i\neq 0$, is a Latin square type srg with parameters
$$(|X|^2, t(|X|-1), |X|+t^2-3t, t^2-t).$$
\end{teor}

Using Theorem~\ref{srg}, one can obtain Latin square type srg from the pseudocyclic association scheme $(\cE, \{\Gamma_a\})$ (the elliptic scheme). These srg have parameters
$$(\frac{1}{2}q^2(q-1)^2, \frac{1}{2}(q-2)(q+1)^2, \frac{1}{2}(3q^2-3q-4), q(q+1)).$$
We note that the Latin square type srg arising from $(\cE, \{\Gamma_a\})$ were mentioned in \cite{tf}, in which another construction of these srg was given.

Now since we have shown that the fusion scheme $(\cE, \{\Delta_a\})$ of the elliptic scheme $(\cE, \{\Gamma_a\})$ is also pseudocyclic when $m$ is an odd prime. We obtain more Latin square type srg via Theorem~\ref{srg}.

\begin{teor}
Let $q=2^m$, where $m$ is an odd prime. Then there exists a Latin square type srg with parameters
$$(\frac{1}{2}q^2(q-1)^2, m(q+1)(\frac{q(q-1)}{2}-1), \lambda, \mu),$$
where $\lambda=\frac{q(q-1)}{2} + m^2(q+1)^2 -3m(q+1)$ and $\mu=m^2(q+1)^2-m(q+1)$.
\end{teor}
\begin{proof} Straightforward.\end{proof}

\noindent{\bf Acknowledgements:} The second author thanks Philips Research Eindhoven, the Netherlands, where part
of this work was carried out. The research of the second author is supported in part by NSF grant DMS 0400411.

\end{document}